\documentclass{article} 
\usepackage{nips11submit_e,times}

\usepackage{graphicx}
\usepackage{amsmath,amssymb} 
\usepackage{algorithm,algorithmic}
\usepackage{color}

\usepackage{epsfig,times}
\usepackage{subfigure}
\usepackage{multirow}

\newtheorem{theorem}{Theorem}

\newtheorem{prop}[theorem]{Proposition}
\newenvironment{proof}[1][Proof]{\begin{trivlist}
\item[\hskip \labelsep {\bfseries #1}]}{\end{trivlist}}

\renewcommand{\dfrac}{\displaystyle \frac} 
\newcommand{\comment}[1]{}

\title{Linearized Alternating Direction Method with Adaptive Penalty for Low-Rank Representation}

\author{
David S.~Hippocampus\thanks{ Use footnote for providing further
information about author (webpage, alternative address)---\emph{not}
for acknowledging
funding agencies.} \\
Department of Computer Science\\
Cranberry-Lemon University\\
Pittsburgh, PA 15213 \\
\texttt{hippo@cs.cranberry-lemon.edu} \\
\And
Coauthor \\
Affiliation \\
Address \\
\texttt{email} \\
\AND
Coauthor \\
Affiliation \\
Address \\
\texttt{email} \\
}

\begin{document}

\maketitle

\begin{abstract}
Low-rank representation (LRR) is an effective method for subspace
clustering and has found wide applications in computer vision and
machine learning. The existing LRR solver is based on the
alternating direction method (ADM). It suffers from $O(n^3)$
computation complexity due to the matrix-matrix multiplications
and matrix inversions, even if partial SVD is used. Moreover,
introducing auxiliary variables also slows down the convergence.
Such a heavy computation load prevents LRR from large scale
applications. In this paper, we generalize ADM by linearizing the
quadratic penalty term and allowing the penalty to change
adaptively. We also propose a novel rule to update the penalty
such that the convergence is fast. With our linearized ADM with
adaptive penalty (LADMAP) method, it is unnecessary to introduce
auxiliary variables and invert matrices. The matrix-matrix
multiplications are further alleviated by using the skinny SVD
representation technique. As a result, we arrive at an algorithm
for LRR with complexity $O(rn^2)$, where $r$ is the rank of the
representation matrix. Numerical experiments verify that for LRR
our LADMAP method is much faster than state-of-the-art algorithms.
Although we only present the results on LRR, LADMAP actually can
be applied to solving more general convex programs.
\end{abstract}

\section{Introduction}
Recently, compressive sensing~\cite{Candes08} and sparse
representation~\cite{Wright-2010-SP} have been hot research topics
and also have found abundant applications in signal processing and
machine learning. Many of the problems in these fields can be
formulated as the following convex programs:
\begin{equation}
\min\limits_{\mathbf{x},\mathbf{y}} f(\mathbf{x}) +
g(\mathbf{y}), \ s.t. \ \mathcal{A}(\mathbf{x}) +
\mathcal{B}(\mathbf{y}) = \mathbf{c},\label{eq:nm}
\end{equation}
where $\mathbf{x}$, $\mathbf{y}$ and $\mathbf{c}$ could be either
vectors or matrices, $f$ and $g$ are convex functions (e.g., the
nuclear norm $\|\cdot\|_*$~\cite{Cai-2008-SVT}, Frobenius norm
$\|\cdot\|$, $l_{2,1}$ norm $\|\cdot\|_{2,1}$~\cite{Liu-2009-L21},
and $l_{1}$ norm $\|\cdot\|_{1}$), and $\mathcal{A}$ and
$\mathcal{B}$ are linear mappings.

As a measure of 2D sparsity, nuclear norm minimization (namely,
$f(\mathbf{X})=\|\mathbf{X}\|_*$) has now attracted a lot of
attention and has been successfully applied to video
processing~\cite{Ji-2010-VideoDenois}, matrix recovery
\cite{Candes-2009-matrix}, unsupervised learning
\cite{Liu-2010-LRR} and semi-supervised
learning~\cite{Goldberg-2010-SemSup}. Typical problems are matrix
completion~\cite{Candes-2009-matrix}, robust principal component
analysis~\cite{Wright-2009-rpca} and their combination
\cite{Candes-2011-RPCA}. A nuclear norm minimization problem could
be reformulated as a semidefinite
program~\cite{Candes-2009-matrix}, hence could be solved by any
off-the-shelf interior point based toolbox, such as CVX. However,
interior point methods cannot handle large scale matrices due to
their $O(n^6)$ complexity in each iteration, where $n\times n$ is
the matrix size. To overcome this issue, several first-order
algorithms have been developed to solve nuclear norm minimization
problems. One method is the singular value thresholding (SVT)
algorithm~\cite{Cai-2008-SVT} which applies the soft-thresholding
operator to the singular values of a certain matrix in each
iteration. The accelerated proximal gradient (APG)
algorithm~\cite{Toh-2009-APG} is also a popular choice due to its
guaranteed $O(k^{-2})$ convergence rate, where $k$ is the
iteration number. The alternating direction method (ADM) has also
regained a lot of attention~\cite{Lin09,Tao11}. It updates the
variables alternately by minimizing the augmented Lagrangian
function with respect to the variables in a Gauss-Seidel manner.

In 2010, Liu et al.~\cite{Liu-2010-LRR} proposed the low-rank
representation (LRR) for robust subspace clustering. Unlike the
sparse representation~\cite{Elhamifar-2009-SSC}, which minimizes
the number of nonzero entries in the representation matrix, LRR
seeks to minimize the rank of the representation matrix. The
mathematical model of LRR is as follows\footnote{Here we switch to
bold capital letters in order to emphasize that the variables are
matrices.}:
\begin{equation}
\min\limits_{\mathbf{Z},\mathbf{E}}\|\mathbf{Z}\|_*+\mu\|\mathbf{E}\|_{2,1},
\ s.t. \ \mathbf{X}=\mathbf{X}\mathbf{Z}+\mathbf{E},\label{eq:lrr}
\end{equation}
where $\mathbf{X}$ is the data matrix. LRR has found wide
applications in computer vision and machine learning, e.g., motion
segmentation, face clustering, and temporal segmentation
\cite{Liu-2010-LRR,Ni-2010-LRRPSD,Favaro-2011-Closed}.

The existing LRR solver~\cite{Liu-2010-LRR} is based on ADM. It
suffers from $O(n^3)$ computation complexity due to the
matrix-matrix multiplications and matrix inversions. Moreover,
introducing auxiliary variables also slows down the convergence,
as there are more variables and constraints. Such a heavy
computation load prevents LRR from large scale applications. In
this paper, we generalize ADM by linearizing the quadratic penalty
term and allowing the penalty to change adaptively. Linearization
makes the auxiliary variables unnecessary, hence waiving the
matrix inversions, while variable penalty makes the convergence
fast. We also propose a novel and simple rule to update the
penalty. We prove the global convergence of linearized ADM with
adaptive penalty (LADMAP) and apply it to LRR, obtaining faster
convergence speed than the original solver. By further
representing $\mathbf{Z}$ as its skinny SVD and utilizing an
advanced functionality of the PROPACK~\cite{Larsen-1998-PROPACK}
package, the complexity of solving LRR by LADMAP becomes only
$O(rn^2)$, as there is no full sized matrix-matrix
multiplications, where $r$ is the rank of the optimal
$\mathbf{Z}$. Although we only present the numerical results on
LRR, LADMAP actually can be applied to solving more general convex
programs.

Our work is inspired by Yang et al.~\cite{Yang-2011-LADM}.
Nonetheless, the difference of our work from theirs is distinct.
First, they only proved the convergence of linearized ADM (LADM)
for a specific problem, namely nuclear norm regularization. Their
proof utilized some special properties of the nuclear norm, while
we prove the convergence of LADM for general problems in
(\ref{eq:nm}). Second, they only proved in the case of fixed
penalty, while we prove in the case of variable penalty. Although
they mentioned the dynamic updating rule proposed in
\cite{He-00-self-adaptive}, their proof cannot be
straightforwardly applied to the case of variable penalty.
Moreover, that rule is for ADM only. Third, the convergence speed
of LADM heavily depends on the choice of penalty. So it is
difficult to choose an optimal fixed penalty that fits for
different problems and problem sizes, while our novel updating
rule for the penalty, although simple, is effective for different
problems and problem sizes.

\section{Linearized Alternating Direction Method with Adaptive Penalty}
\subsection{The Alternating Direction Method}
ADM is now very popular in solving large scale sparse
representation problems~\cite{Boyd-2010-Distributed}. When solving
(\ref{eq:nm}) by ADM, one operates on the following augmented
Lagrangian function:
\begin{equation}
\mathcal{L}_A(\mathbf{x},\mathbf{y},\mathbf{\lambda}) =
f(\mathbf{x}) + g(\mathbf{y}) +
\langle\mathbf{\lambda},\mathcal{A}(\mathbf{x}) +
\mathcal{B}(\mathbf{y})-\mathbf{y}\rangle +
\frac{\beta}{2}\|\mathcal{A}(\mathbf{x}) +
\mathcal{B}(\mathbf{y})-\mathbf{c}\|^2,\label{eq:alf}
\end{equation}
where $\mathbf{\lambda}$ is the Lagrange multiplier,
$\langle\cdot,\cdot\rangle$ is the inner product, and $\beta >0$
is the penalty parameter.
The usual augmented Lagrange multiplier method is to minimize
$\mathcal{L}_A$ w.r.t. $\mathbf{x}$ and $\mathbf{y}$
simultaneously. This is usually difficult and does not exploit the
fact that the objective function is separable. To remedy this
issue, ADM decomposes the minimization of $\mathcal{L}_A$ w.r.t.
$(\mathbf{x},\mathbf{y})$ into two subproblems that minimize
w.r.t. $\mathbf{x}$ and $\mathbf{y}$, respectively. More
specifically, the iterations of ADM go as follows:
\begin{eqnarray}
\mathbf{x}_{k+1}&=&\arg\min\limits_{\mathbf{x}}\mathcal{L}_A(\mathbf{x},\mathbf{y}_k,\mathbf{\lambda}_k),\label{eq:update_x}\\
\mathbf{y}_{k+1}&=&\arg\min\limits_{\mathbf{y}}
\mathcal{L}_A(\mathbf{x}_{k+1},\mathbf{y},\mathbf{\lambda}_k),\label{eq:update_y}\\
\mathbf{\lambda}_{k+1} &=& \mathbf{\lambda}_k +
\beta[\mathcal{A}(\mathbf{x}_{k+1}) +
\mathcal{B}(\mathbf{y}_{k+1})-\mathbf{c}].\label{eq:update_lambda}
\end{eqnarray}
In compressive sensing and sparse representation, as $f$ and $g$
are usually matrix or vector norms, the subproblems
(\ref{eq:update_x}) and (\ref{eq:update_y}) usually have closed
form solutions when $\mathcal{A}$ and $\mathcal{B}$ are identities
\cite{Cai-2008-SVT,Liu-2010-LRR,Yang-2010-ADMl1}. In this case,
ADM is appealing. However, in many problems $\mathcal{A}$ and
$\mathcal{B}$ are not identities. For example, in matrix
completion $\mathcal{A}$ can be a selection matrix, and in 1D
sparse representation $\mathcal{A}$ can be a general matrix. In
this case, there are no closed form solutions to
(\ref{eq:update_x}) and (\ref{eq:update_y}). To overcome this
difficulty, a common strategy is to introduce auxiliary variables
\cite{Liu-2010-LRR,Boyd-2010-Distributed} $\mathbf{u}$ and
$\mathbf{v}$ and reformulate problem (\ref{eq:nm}) into an
equivalent one:
\begin{equation}
\min\limits_{\mathbf{x},\mathbf{y},\mathbf{u},\mathbf{v}} f(\mathbf{x}) +
g(\mathbf{y}), \ s.t. \ \mathcal{A}(\mathbf{u}) +
\mathcal{B}(\mathbf{v}) = \mathbf{c}, \mathbf{x}=\mathbf{u},\mathbf{y}=\mathbf{v},\label{eq:nm'}
\end{equation}
and the corresponding ADM iterations analogous to
(\ref{eq:update_x})-(\ref{eq:update_lambda}) can be deduced. With
more variables and more constraints, the convergence of ADM
becomes slower. Moreover, to update $\mathbf{u}$ and $\mathbf{v}$,
whose subproblems are least squares problems, matrix inversions
are often necessary.

\subsection{Linearized ADM}
To avoid introducing auxiliary variables and still solve
subproblems (\ref{eq:update_x}) and (\ref{eq:update_y})
efficiently, inspired by Yang et al.~\cite{Yang-2011-LADM}, we
propose a linearization technique for (\ref{eq:update_x}) and
(\ref{eq:update_y}). To further accelerate the convergence of the
algorithm, we also propose an adaptive rule for updating the
penalty parameter.

With minor algebra, one can see that subproblem (\ref{eq:update_x}) is equivalent to
\begin{equation}
\mathbf{x}_{k+1}=\arg\min\limits_{\mathbf{x}}f(\mathbf{x})
+ \frac{\beta}{2}\|\mathcal{A}(\mathbf{x}) +
\mathcal{B}(\mathbf{y}_{k}) -
\mathbf{c}+\lambda_k/\beta\|^2.\label{eq:nm_x}
\end{equation}
By linearizing the quadratic term in (\ref{eq:nm_x}) at
$\mathbf{x}_k$ and adding a proximal term, we have the following
approximation:
\begin{equation}
\begin{array}{l}
\mathbf{x}_{k+1}=\arg\min\limits_{\mathbf{x}}f(\mathbf{x})+\langle
\mathcal{A}^*(\lambda_k)+\beta\mathcal{A}^*(\mathcal{A}(\mathbf{x}_k)+\mathcal{B}(\mathbf{y}_k)-\mathbf{c}),
\mathbf{x}-\mathbf{x}_k\rangle
+\frac{\beta\eta_A}{2}\|\mathbf{x}-\mathbf{x}_k\|^2\\
=\arg\min\limits_{\mathbf{x}}f(\mathbf{x})
+\dfrac{\beta\eta_A}{2}\|\mathbf{x}-\mathbf{x}_k+\mathcal{A}^*(\lambda_k+\beta(\mathcal{A}(\mathbf{x}_k)
+ \mathcal{B}(\mathbf{y}_{k})-\mathbf{c}))/(\beta\eta_A)\|^2,
\label{eq:nm_x1}
\end{array}
\end{equation}
where $\mathcal{A}^*$ is the adjoint of $\mathcal{A}$ and
$\eta_A>0$ is a parameter whose proper value will be analyzed
later. The above approximation resembles that of APG
\cite{Toh-2009-APG}, but we do not use APG to solve
(\ref{eq:update_x}) iteratively.

Similarly, subproblem (\ref{eq:update_y}) can be approximated by
\begin{equation}
\mathbf{y}_{k+1}=\arg\min\limits_{\mathbf{y}}g(\mathbf{y})
+\frac{\beta\eta_B}{2}\|\mathbf{y}-\mathbf{y}_k+\mathcal{B}^*(\lambda_k+
\beta(\mathcal{A}(\mathbf{x}_{k+1})+
\mathcal{B}(\mathbf{y}_{k})-\mathbf{c}))/(\beta\eta_B)\|^2.\label{eq:nm_y}
\end{equation}
The update of Lagrange multiplier still goes as
(\ref{eq:update_lambda})\footnote{As in~\cite{Yang-2011-LADM}, we
can also introduce a parameter $\gamma$ and update $\lambda$ as
$\lambda_{k+1} = \lambda_k +
\gamma\beta[\mathcal{A}(\mathbf{x}_{k+1}) +
\mathcal{B}(\mathbf{y}_{k+1})-\mathbf{c}]$. We choose not to do so
in this paper in order not to make the exposition of LADMAP too
complex. The reviewers can refer to Supplementary Material for
full details.}.

\subsection{Adaptive Penalty}
In previous ADM and LADM approaches
\cite{Tao11,Yang-2010-ADMl1,Yang-2011-LADM}, the penalty parameter
$\beta$ is fixed. Some scholars have observed that ADM with a
fixed $\beta$ can converge very slowly and it is nontrivial to
choose an optimal fixed $\beta$. So is LADM. Thus a dynamic
$\beta$ is preferred in real applications. Although Tao et al.
\cite{Tao11} and Yang et al.~\cite{Yang-2011-LADM} mentioned He et
al.'s adaptive updating rule~\cite{He-00-self-adaptive} in their
papers, the rule is for ADM only. We propose the following
adaptive updating strategy for the penalty parameter $\beta$:
\begin{equation}
\beta_{k+1} = \min(\beta_{\max},\rho\beta_k),\label{eq:nm_b}
\end{equation}
where $\beta_{\max}$ is an upper bound of $\{\beta_k\}$. The value of $\rho$ is defined as
\begin{equation}
\rho = \left\{\begin{array}{ll}
\rho_0, & \mbox{if} \ \beta_{k}\max(\sqrt{\eta_A}\|\mathbf{x}_{k+1}-\mathbf{x}_k\|,
\sqrt{\eta_B}\|\mathbf{y}_{k+1}-\mathbf{y}_k\|)/\|\mathbf{c}\|
< \varepsilon_2,\\
1, & \mbox{otherwise},
\end{array}\right.\label{eq:nm_rho}
\end{equation}
where $\rho_0\geq 1$ is a constant. The condition to assign
$\rho=\rho_0$ comes from the analysis on the stopping criteria
(see Section~\ref{sec:stopping}). We recommend that $\beta_0
=\alpha \varepsilon_2$, where $\alpha$ depends on the size of
$\mathbf{c}$. Our updating rule is fundamentally different from He
et al.'s for ADM~\cite{He-00-self-adaptive}, which aims at
balancing the errors in the stopping criteria and involves several
parameters.

\comment{The initial $\beta_0$ is recommended to be chosen as
\begin{equation}
\beta_0=\varepsilon_2 \|\mathbf{c}\|/(2\max(\sqrt{\eta_A}
L_x,\sqrt{\eta_B} L_y)),
\end{equation}
such that the first case of (\ref{eq:nm_rho}) is satisfied, where $L_x$ and $L_y$ are
the estimated upper bound of $\{\|\mathbf{x}_k\|\}$ and $\{\|\mathbf{y}_k\|\}$ (see Proposition~\ref{prop:conv}
for their boundedness), respectively,
which need not be accurate.}

\subsection{Convergence of LADMAP}
To prove the convergence of LADMAP, we first have the following
propositions.
\begin{prop}\label{prop:subgradients}
\begin{equation}
-\beta_k\eta_A(\mathbf{x}_{k+1}-\mathbf{x}_k)-\mathcal{A}^*(\tilde{\lambda}_{k+1}) \in \partial f(\mathbf{x}_{k+1}), \
-\beta_k\eta_B(\mathbf{y}_{k+1}-\mathbf{y}_k)-\mathcal{B}^*(\hat{\lambda}_{k+1})
\in\partial g(\mathbf{y}_{k+1}),\label{eq:kkt_k}
\end{equation}
where
$\tilde{\lambda}_{k+1}=\lambda_k+\beta_k[\mathcal{A}(\mathbf{x}_k)+\mathcal{B}(\mathbf{y}_k)-\mathbf{c}]$,
$\hat{\lambda}_{k+1}=\lambda_k+\beta_k[\mathcal{A}(\mathbf{x}_{k+1})+\mathcal{B}(\mathbf{y}_k)-\mathbf{c}]$,
and $\partial f$ and $\partial g$ are subgradients of $f$ and $g$,
respectively.
\end{prop}
This can be easily proved by checking the optimality conditions of
(\ref{eq:nm_x1}) and (\ref{eq:nm_y}).
\begin{prop}
Denote the operator norms of $\mathcal{A}$ and $\mathcal{B}$ as
$\|\mathcal{A}\|$ and $\|\mathcal{B}\|$, respectively. If
$\{\beta_k\}$ is non-decreasing and upper bounded, $\eta_A >
\|\mathcal{A}\|^2$, $\eta_B > \|\mathcal{B}\|^2$, and
$(x^*,y^*,\lambda^*)$ is any Karush-Kuhn-Tucker (KKT) point of
problem (\ref{eq:nm}) (see (\ref{eq:KKT1})-(\ref{eq:KKT2})), then:
\textbf{(1)}.
$\{\eta_A\|\mathbf{x}_k-\mathbf{x}^*\|^2-\|\mathcal{A}(\mathbf{x}_k-\mathbf{x}^*)\|^2
+\eta_B\|\mathbf{x}_k -
\mathbf{x}^*\|^2+\beta_k^{-2}\|\lambda_{k}-\lambda^*\|^2\}$ is
non-increasing. \textbf{(2)}.\comment{
$\{(\mathbf{x}_k,\mathbf{y}_k,\lambda_k)\}$ is bounded.
\textbf{(3)}.} $\|\mathbf{x}_{k+1}-\mathbf{x}_k\|\to 0$,
$\|\mathbf{y}_{k+1}-\mathbf{y}_k\|\to 0$,
$\|\lambda_{k+1}-\lambda_k\|\to 0$. \comment{\textbf{(4)}.
$\langle \mathbf{x}_{k+1}-\mathbf{x}^*, -
\beta_k\eta_A(\mathbf{x}_{k+1}-\mathbf{x}_k)-\mathcal{A}^*(
\tilde{\lambda}_{k+1})+\mathcal{A}^*(\lambda^*)\rangle$ and
$\langle \mathbf{y}_{k+1}-\mathbf{y}^*, -
\beta_k\eta_B(\mathbf{y}_{k+1}-\mathbf{y}_k)-\mathcal{B}^*(
\hat{\lambda}_{k+1})+\mathcal{B}^*(\lambda^*)\rangle$ are both
nonnegative and convergent to 0.} \label{prop:conv}
\end{prop}
The proof can be found in Supplementary Material. Then we can
prove the convergence of LADMAP, as stated in the following
theorem.
\begin{theorem}
If $\{\beta_k\}$ is non-decreasing and upper bounded, $\eta_A > \|\mathcal{A}\|^2$, and $\eta_B > \|\mathcal{B}\|^2$, then the sequence $\{(\mathbf{x}_k,\mathbf{y}_k,\lambda_k)\}$ generated
by LADMAP converges to a KKT point of problem (\ref{eq:nm}).\label{thm:conv}
\end{theorem}
The proof can be found in Appendix~\ref{sec:proof}.

\subsection{Stopping Criteria}\label{sec:stopping}
\comment{Some existing work, e.g.,
\cite{Liu-2010-LRR,Favaro-2011-Closed}, proposed stopping criteria
out of intuition only, which may not guarantee that the correct
solution is approached. Recently, Lin et al.~\cite{Lin09} and Boyd
et al.~\cite{Boyd-2010-Distributed} suggested that the stopping
criteria be derived from the KKT conditions of a problem. We also
adopt such a strategy.}

The KKT conditions of problem (\ref{eq:nm}) are that there exists
a triple $(\mathbf{x}^*,\mathbf{y}^*,\mathbf{\lambda}^*)$ such
that
\begin{eqnarray}
&\mathcal{A}(\mathbf{x}^*) + \mathcal{B}(\mathbf{y}^*) - \mathbf{c} = \mathbf{0},&\label{eq:KKT1}\\
&-\mathcal{A}^*(\mathbf{\lambda}^*)\in\partial
f(\mathbf{x}^*),-\mathcal{B}^*(\mathbf{\lambda}^*)\in\partial
g(\mathbf{y}^*).&\label{eq:KKT2}
\end{eqnarray}
The triple $(\mathbf{x}^*,\mathbf{y}^*,\mathbf{\lambda}^*)$ is
called a KKT point. So the first stopping criterion is the
feasibility:
\begin{equation}
\|\mathcal{A}(\mathbf{x}_{k+1}) + \mathcal{B}(\mathbf{y}_{k+1}) - \mathbf{c}\|/\|\mathbf{c}\| < \varepsilon_1.\label{eq:stopping1}
\end{equation}

As for the second KKT condition, we rewrite the second part of Proposition~\ref{prop:subgradients} as follows
\begin{equation}
\begin{array}{l}
-\beta_k[\eta_B(\mathbf{y}_{k+1}-\mathbf{y}_k)+\mathcal{B}^*(\mathcal{A}(\mathbf{x}_{k+1}-\mathbf{x}_k))]
-\mathcal{B}^*(\tilde{\lambda}_{k+1})\in\partial
g(\mathbf{y}_{k+1}).
\end{array}
\end{equation}
So for $\tilde{\lambda}_{k+1}$ to satisfy the second KKT
condition, both $\beta_k\eta_A\|\mathbf{x}_{k+1}-\mathbf{x}_k\|$
and
$\beta_k\|\eta_B(\mathbf{y}_{k+1}-\mathbf{y}_k)+\mathcal{B}^*(\mathcal{A}(\mathbf{x}_{k+1}-\mathbf{x}_k))\|$
should be small enough. This leads to the second stopping
criterion:
\begin{equation}
\beta_k\max(\eta_A\|\mathbf{x}_{k+1}-\mathbf{x}_k\|/\|\mathcal{A}^*(\mathbf{c})\|,
\eta_B\|\mathbf{y}_{k+1}-\mathbf{y}_k\|/\|\mathcal{B}^*(\mathbf{c})
\|)\leq\varepsilon_2'.
\end{equation}
By estimating $\|\mathcal{A}^*(\mathbf{c})\|$ and
$\|\mathcal{B}^*(\mathbf{c})\|$ by $\sqrt{\eta_A}\|\mathbf{c}\|$
and $\sqrt{\eta_B}\|\mathbf{c}\|$, respectively, we arrive at the
second stopping criterion in use:
\begin{equation}
\beta_k\max(\sqrt{\eta_A}\|\mathbf{x}_{k+1}-\mathbf{x}_k\|,
\sqrt{\eta_B}\|\mathbf{y}_{k+1}-\mathbf{y}_k\|)/\|\mathbf{c}\|\leq\varepsilon_2.\label{eq:stopping2}
\end{equation}

Finally, we summarize our LADMAP algorithm in
Algorithm~\ref{alg:ialm_nnm}.
\begin{algorithm}[tb]
   \caption{LADMAP for Problem (\ref{eq:nm})}
   \label{alg:ialm_nnm}
\begin{algorithmic}
   \STATE {\bfseries Initialize:}
     Set $\varepsilon_1>0$, $\varepsilon_2 > 0$, $\beta_{\max}\gg\beta_0 > 0$,
   $\eta_A > \|\mathcal{A}\|^2$, $\eta_B > \|\mathcal{B}\|^2$, $\mathbf{x}_0$, $\mathbf{y}_0$, $\mathbf{\lambda}_0$, and $k\leftarrow 0$.
   \WHILE {(\ref{eq:stopping1}) or (\ref{eq:stopping2}) is not satisfied}
   \STATE {\bfseries Step 1:} Update $\mathbf{x}$ by solving
   (\ref{eq:nm_x1}).
   \STATE {\bfseries Step 2:} Update $\mathbf{y}$ by solving
   (\ref{eq:nm_y}).
   \STATE {\bfseries Step 3:} Update $\lambda$ by (\ref{eq:update_lambda}).
   \STATE {\bfseries Step 4:} Update $\beta$ by (\ref{eq:nm_b}) and (\ref{eq:nm_rho}).
   \STATE {\bfseries Step 5:} $k\leftarrow k+1$.
   \ENDWHILE
\end{algorithmic}
\end{algorithm}

\section{Applying LADMAP to LRR}
In this section, we apply LADMAP to solving the LRR problem
(\ref{eq:lrr}). We also compare LADMAP with other state-of-the-art
algorithms for LRR. The reason we choose LRR as an example of
applications of LADMAP is twofold. First, LRR has become an
important mathematical model in machine learning. Second, unlike
other established nuclear norm minimization problems, such as
matrix completion~\cite{Candes-2009-matrix} and robust principal
component analysis~\cite{Candes-2011-RPCA}, if not carefully
treated, the complexity of solving LRR is still $O(n^3)$, even if
partial SVD is used.

\subsection{Solving LRR by LADMAP}\label{sec:lrr_ladmap}
As the LRR problem (\ref{eq:lrr}) is a special case of problem
(\ref{eq:nm}), LADMAP can be directly applied to it. The two
subproblems both have closed form solutions. In the subproblem for
updating $\mathbf{E}$, one may apply the $l_{2,1}$-norm shrinkage
operator~\cite{Liu-2010-LRR}, with a threshold $\beta_k^{-1}$, to
matrix $\mathbf{M}_k=-\mathbf{X}\mathbf{Z}_k
+\mathbf{X}-\mathbf{\Lambda}_k/\beta_k$. In the subproblem for
updating $\mathbf{Z}$, one has to apply the singular value
shrinkage operator~\cite{Cai-2008-SVT}, with a threshold
$(\beta_k\eta_X)^{-1}$, to matrix
$\mathbf{N}_k=\mathbf{Z}_k-\eta_X^{-1}
\mathbf{X}^T(\mathbf{X}\mathbf{Z}_k
+\mathbf{E}_{k+1}-\mathbf{X}+\mathbf{\Lambda}_k/\beta_k)$, where
$\eta_X
> \sigma_{\max}^2(\mathbf{X})$. If $\mathbf{N}_k$ is formed
explicitly, the usual technique of partial SVD, using PROPACK
\cite{Larsen-1998-PROPACK}, with rank prediction
\cite{Lin09,Yang-2011-LADM,Toh-2009-APG} can be utilized to
compute the leading $r$ singular values and associated vectors of
$\mathbf{N}_k$ efficiently, making the complexity of SVD
computation $O(rn^2)$, where $r$ is the predicted rank of
$\mathbf{Z}_{k+1}$ and $n$ is the column number of $\mathbf{X}$.
Note that as $\beta_k$ is non-decreasing, the predicted rank is
almost non-decreasing, making the iterations computationally
efficient.

Up to now, LADMAP for LRR is still of complexity $O(n^3)$,
although partial SVD is already used. This is because forming
$\mathbf{M}_k$ and $\mathbf{N}_k$ requires full sized
matrix-matrix multiplications, e.g., $\mathbf{X}\mathbf{Z}_k$. To
break this complexity bound, we introduce another technique to
further accelerate LADMAP for LRR. By representing $\mathbf{Z}_k$
as its skinny SVD: $\mathbf{Z}_k =
\mathbf{U}_k\mathbf{\Sigma}_k\mathbf{V}_k^T$, some of the full
sized matrix-matrix multiplications are gone: they are replaced by
successive reduced sized matrix-matrix multiplications. For
example, when updating $\mathbf{E}$, $\mathbf{X}\mathbf{Z}_k$ is
computed as
$((\mathbf{X}\mathbf{U}_k)\mathbf{\Sigma}_k)\mathbf{V}_k^T$,
reducing the complexity to $O(rn^2)$. When computing the partial
SVD of $\mathbf{N}_k$, things are more complicated. If we form
$\mathbf{N}_k$ explicitly, we will face with computing
$\mathbf{X}^T(\mathbf{X}+\mathbf{\Lambda}_k/\beta_k)$, which is
neither low-rank nor sparse\footnote{When forming $\mathbf{N}_k$
explicitly, $\mathbf{X}^T\mathbf{X}\mathbf{Z}_k$ can be computed
as
$((\mathbf{X}^T(\mathbf{X}\mathbf{U}_k))\mathbf{\Sigma}_k)\mathbf{V}_k^T$,
whose complexity is still $O(rn^2)$, while
$\mathbf{X}^T\mathbf{E}_{k+1}$ could also be accelerated as
$\mathbf{E}_{k+1}$ is a column-sparse matrix.}. Fortunately, in
PROPACK the bi-diagonalizing process of $\mathbf{N}_k$ is done by
the Lanczos procedure~\cite{Larsen-1998-PROPACK}, which only
requires to compute matrix-vector multiplications
$\mathbf{N}_k\mathbf{v}$ and $\mathbf{u}^T\mathbf{N}_k$, where
$\mathbf{u}$ and $\mathbf{v}$ are some vectors in the Lanczos
procedure. So we may compute $\mathbf{N}_k\mathbf{v}$ and
$\mathbf{u}^T\mathbf{N}_k$ by multiplying the vectors $\mathbf{u}$
and $\mathbf{v}$ successively with the component matrices in
$\mathbf{N}_k$, rather than forming $\mathbf{N}_k$ explicitly. So
the computation complexity of partial SVD of $\mathbf{N}_k$ is
still $O(rn^2)$. Consequently, with our acceleration techniques,
the complexity of our accelerated LADMAP for LRR is $O(rn^2)$. The
accelerated LADMAP is summarized in
Algorithm~\ref{alg:LADMAP_lrr}.

\begin{algorithm}[tb]
   \caption{Accelerated LADMAP for LRR (\ref{eq:lrr})}
   \label{alg:LADMAP_lrr}
\begin{algorithmic}
   \STATE {\bfseries Input:} Observation matrix $\mathbf{X}$ and parameter $\mu >0$.
   \STATE {\bfseries Initialize:} Set $\mathbf{E}_0$, $\mathbf{Z}_0$ and $\mathbf{\Lambda}_0$
   to zero matrices, where $\mathbf{Z}_0$ is represented as $(\mathbf{U}_0,\mathbf{\Sigma}_0,\mathbf{V}_0)\leftarrow(\mathbf{0},\mathbf{0},\mathbf{0})$.
   Set $\varepsilon_1>0$, $\varepsilon_2>0$, $\beta_{\max}\gg\beta_0>0$, $\eta_X>\sigma_{\max}^2(\mathbf{X})$, and $k\leftarrow 0$.
   \WHILE {(\ref{eq:stopping1}) or (\ref{eq:stopping2}) is not satisfied}
   \STATE {\bfseries Step 1:} Update $\mathbf{E}_{k+1}=\arg\min\limits_{\mathbf{E}}\mu\|\mathbf{E}\|_{2,1}+
   \frac{\beta_k}{2}\|\mathbf{E}+(\mathbf{X}\mathbf{U}_k)\mathbf{\Sigma}_k\mathbf{V}_k^T-\mathbf{X}+\mathbf{\Lambda}_k/\beta_k\|^2$.
   This subproblem can be solved by using Lemma~3.2 in~\cite{Liu-2010-LRR}.
   \STATE {\bfseries Step 2:} Update the skinny SVD $(\mathbf{U}_{k+1},\mathbf{\Sigma}_{k+1},\mathbf{V}_{k+1})$ of
   $\mathbf{Z}_{k+1}$. First, compute the partial SVD
   $\tilde{\mathbf{U}}_r\tilde{\mathbf{\Sigma}}_r\tilde{\mathbf{V}}_r^T$ of the \emph{implicit} matrix $\mathbf{N}_k$,
   which is bi-diagonalized by the successive matrix-vector multiplication technique
   described in Section~\ref{sec:lrr_ladmap}, and the rank $r$ is predicted as in~\cite{Lin09,Yang-2011-LADM,Toh-2009-APG}.
   Second,
   $\mathbf{U}_{k+1}= \tilde{\mathbf{U}}_r(:,1:r')$, $\mathbf{\Sigma}_{k+1}=
   \tilde{\mathbf{\Sigma}}_r(1:r',1:r')-(\beta_k\eta_X)^{-1}\mathbf{I}$,
   $\mathbf{V}_{k+1}= \tilde{\mathbf{V}}_r(:,1:r')$, where $r'$ is the
   number of singular values in $\mathbf{\Sigma}_r$ that are greater than
   $(\beta_k\eta_X)^{-1}$.
   \STATE {\bfseries Step 3:} Update
   $\mathbf{\Lambda}_{k+1}=\mathbf{\Lambda}_k+\beta_{k}((\mathbf{X}\mathbf{U}_{k+1})\mathbf{\Sigma}_{k+1}\mathbf{V}_{k+1}^T+\mathbf{E}_{k+1}-\mathbf{X})$.
   \STATE {\bfseries Step 4:} Update $\beta_{k+1}$ by (\ref{eq:nm_b})-(\ref{eq:nm_rho}).
   \STATE {\bfseries Step 5:} $k\leftarrow k+1$.
   \ENDWHILE
\end{algorithmic}
\end{algorithm}

\subsection{Comparison with Other Methods}
As shown in~\cite{Liu-2010-LRR}, the LRR problem can be solved by
the classic ADM. However, their algorithm requires an auxiliary
variable, matrix- matrix multiplication and inversion of matrices,
resulting in $O(n^3)$ computation complexity and slow convergence.

The LRR problem can also be solved \emph{approximately} by being
reformulated to the following unconstrained optimization problem:
$
\min\limits_{\mathbf{Z},\mathbf{E}}\beta(\|\mathbf{Z}\|_*+\mu\|\mathbf{E}\|_{2,1})
+\frac{1}{2}\|\mathbf{X}-\mathbf{X}\mathbf{R}-\mathbf{E}\|^2 $,
where $\beta > 0$ is a relaxation parameter. Then APG with the
continuation technique~\cite{Toh-2009-APG}, which is to reduce
$\beta$ gradually by
$\beta_{k+1}=\max(\beta_{\min},\theta\beta_k)$, can be applied to
solve this problem\footnote{Please see Supplementary Material for
the detail of solving LRR by APG.}. \comment{A continuation
technique~\cite{Toh-2009-APG}, which varies $\beta$ by starting
from a large initial value and decreasing it geometrically until
reaching a lower bound $\beta_{\min}$, can significantly speed up
the convergence. }Compared with APG, which can only find an
approximate solution, LADMAP can produce a much more accurate
solution as it is proven to converge to an exact solution.

The linearization technique has also been applied to other
optimization methods. For example, Yin~\cite{Yin-2010-LSB} applied
this technique to the Bregman iteration for solving compressive
sensing problems and proved that the linearized Bregman method
converges to an exact solution \emph{conditionally}. In
comparison, LADMAP always converges to an exact solution.

\section{Experimental Results}


In this section, we report numerical results on the standard
LADMAP, the accelerated LADMAP and other state-of-the-art
algorithms, including APG, ADM\footnote{We use the Matlab code
provided online by the authors of~\cite{Liu-2010-LRR}.} and LADM,
for LRR based data clustering problems. APG, ADM, LADM and LADMAP
all utilize the Matlab version of PROPACK
\cite{Larsen-1998-PROPACK}. For the accelerated LADMAP, we provide
two function handles to PROPACK which fulfils the successive
matrix-vector multiplications. All experiments are run and timed
on a PC with an Intel Core i5 CPU at 2.67GHz and with 4GB of
memory, running Windows 7 and Matlab version 7.10.

We test and compare these solvers on both synthetic multiple
subspaces data and the real world motion data (Hopkin155 motion
segmentation database~\cite{Tron-2007-Hop}). For APG, we set the
parameters $\beta_0 = 0.01$, $\beta_{\min}=10^{-10}$, $\theta =
0.9$ and the Lipschitz constant $\tau =
\sigma_{\max}^2(\mathbf{X})$. The parameters of ADM and LADM are
the same as that in~\cite{Liu-2010-LRR} and~\cite{Yang-2011-LADM},
respectively. In particular, for LADM the penalty is fixed at
$\beta=2.5/\min(m,n)$, where $m\times n$ is the size of
$\mathbf{X}$. For LADMAP, we set $\varepsilon_1 = 10^{-4}$,
$\varepsilon_2=10^{-5}$, $\beta_0 = \min(m,n)\varepsilon_2$,
$\beta_{\max}=10^{10}$, $\rho_0 = 1.9$, and $\eta_X = 1.02
\sigma_{\max}^2(\mathbf{X})$. As the code of ADM was downloaded,
its stopping criteria,
$\|\mathbf{X}\mathbf{Z}_k+\mathbf{E}_k-\mathbf{X}\|/\|\mathbf{X}\|\leq\varepsilon_1$
and
$\max(\|\mathbf{E}_{k}-\mathbf{E}_{k-1}\|/\|\mathbf{X}\|,\|\mathbf{Z}_{k}-\mathbf{Z}_{k-1}\|/\|\mathbf{X}\|)\leq\varepsilon_2$,
are used in all our experiments\footnote{Note that the second
criterion differs from that in (\ref{eq:stopping2}). However, this
does not harm the convergence of LADMAP because
(\ref{eq:stopping2}) is always checked when updating $\beta_{k+1}$
(see (\ref{eq:nm_rho})).}.

\comment{As the code of ADM was downloaded, its stopping criteria,
$\mbox{FeasError}\leq\varepsilon_1$ and
$\max(\mbox{RelChg}(\mathbf{E}_k),\mbox{RelChg}(\mathbf{Z}_k))\leq\varepsilon_2$,
are used in all our experiments, where the feasibility error
(FeasError) and the relative changes (RelChg) are defined as
$\mbox{FeasError}=\|\mathbf{X}\mathbf{Z}_k+\mathbf{E}_k-\mathbf{X}\|/\|\mathbf{X}\|$
and
$\mbox{RelChg}(\mathbf{E}_k)=\|\mathbf{E}_{k}-\mathbf{E}_{k-1}\|/\|\mathbf{X}\|$,
$\mbox{RelChg}(\mathbf{Z}_k)=\|\mathbf{Z}_{k}-\mathbf{Z}_{k-1}\|/\|\mathbf{X}\|$,
respectively\footnote{Note that the relative change criteria
differ from that in (\ref{eq:stopping2}). However, this does not
harm the convergence of LADMAP because (\ref{eq:stopping2}) is
always checked when updating $\beta_{k+1}$ (see
(\ref{eq:nm_rho})).}.}

\subsection{On Synthetic Data}
The synthetic test data, parameterized as ($s$, $p$, $d$,
$\tilde{r}$), is created by the same procedure in
\cite{Liu-2010-LRR}. $s$ independent subspaces
$\{\mathcal{S}_i\}_{i=1}^s$ are constructed, whose bases
$\{\mathbf{U}_i\}_{i=1}^s$ are generated by
$\mathbf{U}_{i+1}=\mathbf{T}\mathbf{U}_i, \ 1 \leq i \leq s-1$,
where $\mathbf{T}$ is a random rotation and $\mathbf{U}_1$ is a $d
\times \tilde{r}$ random orthogonal matrix. So each subspace has a
rank of $\tilde{r}$ and the data has an ambient dimension of $d$.
Then $p$ data points are sampled from each subspace by
$\mathbf{X}_i=\mathbf{U}_i\mathbf{Q}_i, \ 1\leq i\leq s$, with
$\mathbf{Q}_i$ being an $\tilde{r} \times p$ i.i.d. zero mean unit
variance Gaussian matrix $\mathcal{N}(0,1)$. 20$\%$ samples are
randomly chosen to be corrupted by adding Gaussian noise with zero
mean and standard deviation $0.1\|\mathbf{x}\|$. We empirically
find that LRR achieves the best clustering performance on this
data set when $\mu=0.1$. So we test all algorithms with $\mu =
0.1$ in this experiment. To measure the relative errors in the
solutions, we run the standard LADMAP 2000 iterations with
$\beta_{\max}=10^{3}$ to establish the ground truth solution
$(\mathbf{E}_0,\mathbf{Z}_0)$.

The computational comparison is summarized in Table~\ref{tab:toy}.
We can see that the iteration numbers and the CPU times of both
the standard and accelerated LADMAP are much less than those of
other methods, and the accelerated LADMAP is further much faster
than the standard LADMAP. Moreover, the advantage of the
accelerated LADMAP is even greater when the ratio $\tilde{r}/p$,
which is roughly the ratio of the rank of $\mathbf{Z}_0$ to the
size of $\mathbf{Z}_0$, is smaller, which testifies to the
complexity estimations on the standard and accelerated LADMAP for
LRR. It is noteworthy that the iteration numbers of ADM and LADM
seem to grow with the problem sizes, while that of LADMAP is
rather constant. Moreover, LADM is not faster than ADM. In
particular, on the last data we were unable to wait until LADM
stopped. Finally, as APG converges to an approximate solution to
(\ref{eq:lrr}), its relative errors are larger and its clustering
accuracy is lower than ADM and LADM based methods. \comment{In
Figure~\ref{fig:toy}, we plot the relative changes of
$\mathbf{E}_k$ and $\mathbf{Z}_k$ and the feasibility errors at
all iterations for four test algorithms, respectively. We can see
the errors of LADMAP in the two KKT conditions drop much quicker
than other methods.}

\begin{table*}[th]
\begin{center}
\caption{Comparison among APG, ADM, LADM, standard LADMAP and
accelerated LADMAP (denoted as LADMAP(A)) on the synthetic data.
For each quadruple ($s$, $p$, $d$, $\tilde{r}$), the LRR problem,
with $\mu = 0.1$, was solved for the same data using different
algorithms. We present typical running time (in $\times10^3$
seconds), iteration number, relative error (\%) of output solution
$(\hat{\mathbf{E}},\hat{\mathbf{Z}})$ and the clustering accuracy
(\%) of tested algorithms, respectively.}\label{tab:toy}
\begin{tabular}{c|c|c|c|c|c|c}\hline
Size ($s$, $p$, $d$, $\tilde{r}$) & Method & Time &  Iter. &
$\frac{\|\hat{\mathbf{Z}}-\mathbf{Z}_0\|}{\|\mathbf{Z}_0\|}$ &
$\frac{\|\hat{\mathbf{E}}-\mathbf{E}_0\|}{\|\mathbf{E}_0\|}$ &
Acc.
\\\hline\hline
\multirow{5}{*}{(10, 20,200, 5)}
& APG       & 0.0332         & 110          & 2.2079          & 1.5096          & 81.5 \\
& ADM       & 0.0529         & 176          & 0.5491          & 0.5093          & \textbf{90.0}\\
& LADM      & 0.0603         & 194          & \textbf{0.5480} & \textbf{0.5024}          & \textbf{90.0}\\
& LADMAP    & 0.0145         & \textbf{46}  & \textbf{0.5480} & \textbf{0.5024} & \textbf{90.0 }\\
& LADMAP(A) & \textbf{0.0010}& \textbf{46}  & \textbf{0.5480} &
\textbf{0.5024} & \textbf{90.0}\\\hline \multirow{5}{*}{(15, 20,300,
5)}
& APG  & 0.0869               & 106         & 2.4824          & 1.0341 & 80.0 \\
& ADM  & 0.1526               & 185         & 0.6519          & 0.4078 & 83.7\\
& LADM & 0.2943               & 363         & \textbf{0.6518} & \textbf{0.4076} & \textbf{86.7} \\
& LADMAP & 0.0336             & \textbf{41} & \textbf{0.6518} & \textbf{0.4076} & \textbf{86.7} \\
& LADMAP(A) & \textbf{0.0015} & \textbf{41} & \textbf{0.6518} &
\textbf{0.4076} & \textbf{86.7} \\\hline \multirow{5}{*}{(20, 25,
500, 5)}
& APG  & 1.8837         & 117          & 2.8905          & 2.4017 & 72.4 \\
& ADM  & 3.7139         & 225          & 1.1191          & 1.0170 & 80.0\\
& LADM & 8.1574         & 508          & \textbf{0.6379} & \textbf{0.4268}& 80.0\\
& LADMAP &  0.7762      & \textbf{40}  & \textbf{0.6379}  & \textbf{0.4268} & \textbf{84.6} \\
& LADMAP(A) & \textbf{0.0053}& \textbf{40} & \textbf{0.6379} &
\textbf{0.4268} & \textbf{84.6}\\\hline \multirow{5}{*}{(30, 30,
900, 5)}
& APG  & 6.1252         & 116          & 3.0667          & 0.9199 & 69.4 \\
& ADM  & 11.7185         & 220          & 0.6865          & 0.4866 & \textbf{76.0}\\
& LADM & N.A.         & N.A.          & N.A. & N.A. & N.A.\\
& LADMAP &  2.3891      & \textbf{44}  & \textbf{0.6864}  & \textbf{0.4294} & \textbf{80.1} \\
& LADMAP(A) & \textbf{0.0058}& \textbf{44} & \textbf{0.6864} &
\textbf{0.4294} & \textbf{80.1}\\\hline
\end{tabular}
\end{center}\vspace{-0.5em}
\end{table*}

\comment{
\begin{figure*}[th]
\begin{center}
\begin{tabular}{c@{\extracolsep{0.2em}}c@{\extracolsep{0.2em}}c}
\includegraphics[width=0.33\textwidth,
keepaspectratio]{figure/5x20x5x100_relChg_E}
&\includegraphics[width=0.33\textwidth,
keepaspectratio]{figure/5x20x5x100_relChg_R}
&\includegraphics[width=0.33\textwidth,
keepaspectratio]{figure/5x20x5x100_recErr}\\
(a) Relative Changes of $\mathbf{E}_k$ & (b) Relative Changes of $\mathbf{Z}_k$ & (c) Feasibility Errors \\
\end{tabular}
\end{center}
\caption{Convergence behaviors of APG, ADM, LADM, LADMAP on the
toy data $\mathbf{X}$ generated with parameters (5, 20, 100, 5).
The changes and errors are in $\log_{10}$ scale. In (a) and (b),
as the relative changes of $\mathbf{E}_k$ and $\mathbf{Z}_k$ in
the first several iterations are zeros, which corresponds to
$-\infty$ in the plots, we only report the nonzero relative
changes of $\mathbf{E}_k$ and $\mathbf{Z}_k$.}\label{fig:toy}
\end{figure*}
}

\begin{table*}[th]
\begin{center}
\caption{Comparison among APG, ADM, LADM, standard LADMAP and
accelerated LADMAP on the Hopkins155 database. We present their
average computing time (in seconds), average number of iterations
and average classification errors (\%) on all 156
sequences.}\label{tab:motion}
\begin{tabular}{c|c|c|c|c|c|c|c|c|c}\hline
&\multicolumn{3}{c|}{Two Motion} & \multicolumn{3}{|c|}{Three
Motion} & \multicolumn{3}{|c}{All}\\\hline & Time & Iter. &  CErr.
& Time & Iter. &  CErr. & Time & Iter. & CErr.\\\hline
APG        & 15.7836        & 90          & 5.77            & 46.4970          & 90         & \textbf{16.52}  & 22.6277         &  90        & 8.36\\
ADM        & 53.3470        & 281         & \textbf{5.72}   & 159.8644         & 284        & \textbf{16.52}  & 77.0864         &  282       & \textbf{8.33} \\
LADM       & 9.6701         & 110         & 5.77            & 22.1467          & 64         & \textbf{16.52}  & 12.4520         &  99        & 8.36\\
LADMAP     & 3.6964          & \textbf{22}     & \textbf{5.72}   & 10.9438          & \textbf{22}& \textbf{16.52}  & 5.3114         & \textbf{22}& \textbf{8.33}\\
LADMAP(A)  & \textbf{2.1348} & \textbf{22}     & \textbf{5.72}   &
\textbf{6.1098} & \textbf{22}         & \textbf{16.52}  &
\textbf{3.0202} & \textbf{22} & \textbf{8.33}\\\hline
\end{tabular}
\end{center}\vspace{-0.5em}
\end{table*}

\subsection{On Real World Data}
We further test the performance of these algorithms on the
Hopkins155 database~\cite{Tron-2007-Hop}. This database consists
of 156 sequences, each of which has 39 to 550 data vectors drawn
from two or three motions. For computational efficiency, we
preprocess the data by projecting it to be 5-dimensional using
PCA. As $\mu=2.4$ is the best parameter for this database
\cite{Liu-2010-LRR}, we test all algorithms with $\mu = 2.4$.

Table~\ref{tab:motion} shows the comparison among APG, ADM, LADM,
standard LADMAP and accelerated LADMAP on this database. We can
also see that the standard and accelerated LADMAP are much faster
than APG, ADM, and LADM, and the accelerated LADMAP is also faster
than the standard LADMAP. However, in this experiment the
advantage of the accelerated LADMAP over the standard LADMAP is
not as dramatic as that in Table~\ref{tab:toy}. This is because on
this data $\mu$ is chosen as $2.4$, which cannot make the rank of
the ground truth solution $\mathbf{Z}_0$ much smaller than the
size of $\mathbf{Z}_0$.

\section{Conclusions}
In this paper, we propose a linearized alternating direction
method with adaptive penalty (LADMAP) and apply it to solving the
LRR problem. With linearization, auxiliary variables need not be
introduced for closed-form solutions, when the objective functions
are matrix or vector norms. Moreover, with fewer variables and
constraints, the convergence also becomes faster. Allowing the
penalty to change adaptively further accelerates the convergence
of LADM. When applying LADMAP to LRR, by representing the
representation matrix as its skinny SVD, full sized matrix-matrix
multiplications are avoided by using successive reduced sized
matrix-matrix multiplications instead, and successive
matrix-vector multiplications are introduced to compute the
partial SVD. Finally, we are able to solve LRR at a computation
complexity of $O(rn^2)$, which is highly advantageous over the
existing LRR solvers. Numerical results demonstrate that LADMAP
converges faster than LADM and ADM and our acceleration techniques
are effective on LRR. Although we only present results on LRR,
LADMAP is actually a general method that can be applied to other
convex programs. We will test it with more problems in sparse
representation in the future.

\appendix
\section{Proof of Theorem~\ref{thm:conv}}\label{sec:proof}
\begin{proof}
By Proposition~\ref{prop:conv} (1),
$\{(\mathbf{x}_k,\mathbf{y}_k,\lambda_k)\}$ is bounded, hence has
an accumulation point, say $
(\mathbf{x}_{k_j},\mathbf{y}_{k_j},\lambda_{k_j}) \to
(\mathbf{x}^{\infty},\mathbf{y}^{\infty},\lambda^{\infty}) $. We
accomplish the proof in two steps.

\textbf{1.} We first prove that
$(\mathbf{x}^{\infty},\mathbf{y}^{\infty},\lambda^{\infty})$ is a
KKT point of problem (\ref{eq:nm}).

By Proposition~\ref{prop:conv} (2), $
\mathcal{A}(\mathbf{x}_{k+1})+\mathcal{B}(\mathbf{y}_{k+1})-\mathbf{c}
= \beta_k^{-1}(\lambda_{k+1}-\lambda_k) \to 0.$ This shows that
any accumulation point of $\{(\mathbf{x}_{k},\mathbf{y}_{k})\}$ is
a feasible solution.

By letting $k=k_j - 1$ in Proposition~\ref{prop:subgradients} and
the definition of subgradient, we have
\begin{eqnarray*}
\begin{array}{rl}
f(\mathbf{x}_{k_j}) + g(\mathbf{y}_{k_j})\leq & f(\mathbf{x}^*) +
g(\mathbf{y}^*)
+\langle \mathbf{x}_{k_j} - \mathbf{x}^*,-\beta_{k_j-1}\eta_A(\mathbf{x}_{k_j}-\mathbf{x}_{k_j-1})-\mathcal{A}^*(\tilde{\lambda}_{k_j}) \rangle \\
&+\langle \mathbf{y}_{k_j} - \mathbf{y}^*,-\beta_{k_j-1}\eta_B(\mathbf{y}_{k_j}-\mathbf{y}_{k_j-1})-\mathcal{B}^*(\hat{\lambda}_{k_j}) \rangle.
\end{array}
\end{eqnarray*}
Let $j\to +\infty$, by observing Proposition~\ref{prop:conv} (2),
we have
\begin{eqnarray*}
\begin{array}{rl}
f(\mathbf{x}^\infty) + g(\mathbf{y}^\infty)&\leq  f(\mathbf{x}^*) + g(\mathbf{y}^*)
+\langle \mathbf{x}^\infty - \mathbf{x}^*,-\mathcal{A}^*(\lambda^\infty) \rangle
+\langle \mathbf{y}^\infty - \mathbf{y}^*,-\mathcal{B}^*(\lambda^\infty) \rangle\\
&=f(\mathbf{x}^*) + g(\mathbf{y}^*)
-\langle \mathcal{A}(\mathbf{x}^\infty - \mathbf{x}^*),\lambda^\infty \rangle
-\langle \mathcal{B}(\mathbf{y}^\infty - \mathbf{y}^*),\lambda^\infty) \rangle\\
&=f(\mathbf{x}^*) + g(\mathbf{y}^*)
-\langle \mathcal{A}(\mathbf{x}^\infty)+\mathcal{B}(\mathbf{y}^\infty) - \mathcal{A}(\mathbf{x}^*)-\mathcal{B}(\mathbf{y}^*),\lambda^\infty \rangle\\
&=f(\mathbf{x}^*) + g(\mathbf{y}^*),
\end{array}
\end{eqnarray*}
where we have used the fact that both $(\mathbf{x}^\infty,\mathbf{y}^\infty)$ and $(\mathbf{x}^*,\mathbf{y}^*)$ are feasible solutions.
So we conclude that
$(\mathbf{x}^{\infty},\mathbf{y}^{\infty})$ is an optimal
solution to (\ref{eq:nm}).

Again, let $k = k_j-1$
in Proposition~\ref{prop:subgradients} and by the definition of subgradient, we have
\begin{equation}
f(\mathbf{x})\geq f(\mathbf{x}_{k_j})
+\langle\mathbf{x}-\mathbf{x}_{k_j},-\beta_{k_j-1}\eta_A(\mathbf{x}_{k_j}-\mathbf{x}_{k_j-1})
-\mathcal{A}^*(\tilde{\lambda}_{k_j})\rangle, \ \forall
\mathbf{x}.
\end{equation}
Fix $\mathbf{x}$ and let $j\to +\infty$, we see that $$
f(\mathbf{x}) \geq f(\mathbf{x}^{\infty}) + \langle
\mathbf{x}-\mathbf{x}^{\infty},
-\mathcal{A}^*(\lambda^{\infty})\rangle, \ \forall\mathbf{x}.$$ So
$-\mathcal{A}^*(\lambda^{\infty})\in\partial
f(\mathbf{x}^{\infty})$. Similarly,
$-\mathcal{B}^*(\lambda^{\infty})\in\partial
g(\mathbf{y}^{\infty})$. Therefore,
$(\mathbf{x}^{\infty},\mathbf{y}^{\infty},\lambda^{\infty})$ is a
KKT point of problem (\ref{eq:nm}).

\textbf{2.} We next prove that the whole sequence
$\{(\mathbf{x}_k,\mathbf{y}_k,\lambda_k)\}$ converges to
$(\mathbf{x}^{\infty},\mathbf{y}^{\infty},\lambda^{\infty})$.

By choosing
$(\mathbf{x}^*,\mathbf{y}^*,\lambda^*)=(\mathbf{x}^{\infty},\mathbf{y}^{\infty},\lambda^{\infty})$
in Proposition~\ref{prop:conv}, we have $
\eta_A\|\mathbf{x}_{k_j}-\mathbf{x}^{\infty}\|^2-\|\mathcal{A}(\mathbf{x}_{k_j}-\mathbf{x}^{\infty})\|^2
+\eta_B\|\mathbf{y}_{k_j} -
\mathbf{y}^{\infty}\|^2+\beta_{k_j}^{-2}\|\lambda_{k_j}-\lambda^{\infty}\|^2
\to 0 $. By Proposition~\ref{prop:conv} (1), we readily have
$\eta_A\|\mathbf{x}_k-\mathbf{x}^{\infty}\|^2-\|\mathcal{A}(\mathbf{x}_k-\mathbf{x}^{\infty})\|^2
+\eta_B\|\mathbf{y}_k -
\mathbf{y}^{\infty}\|^2+\beta_k^{-2}\|\lambda_{k}-\lambda^{\infty}\|^2
\to 0 $. So $(\mathbf{x}_k,\mathbf{y}_k,\lambda_k) \to
(\mathbf{x}^{\infty},\mathbf{y}^{\infty},\lambda^{\infty})$.

As $(\mathbf{x}^{\infty},\mathbf{y}^{\infty},\lambda^{\infty})$
can be an arbitrary accumulation point of
$\{(\mathbf{x}_k,\mathbf{y}_k,\lambda_k)\}$, we may conclude that
$\{(\mathbf{x}_k,\mathbf{y}_k,\lambda_k)\}$ converges to a KKT
point of problem (\ref{eq:nm}).
\end{proof}

\bibliographystyle{plain}
\bibliography{ladm}

\end{document}